\documentclass [11pt,a4paper]{article}
\usepackage{amssymb,amscd}
\usepackage{amsmath}
\usepackage{pifont}
\usepackage{stmaryrd}
\usepackage{amsfonts}
\usepackage{bbm}
\usepackage{amssymb}
\usepackage{lipsum}
\usepackage{ntheorem}
\usepackage[affil-it]{authblk}
 \theoremstyle{plain}
\newtheorem{theo}{Theorem}[section]
\newcounter{tmp}
\newtheorem{corollary}[theo]{Corollary}

\newtheorem*{corollary-non}{Corollary}[section]

\newtheorem*{Cj-non}{Conjecutre}[section]

\newtheorem{lemma}[theo]{Lemma}

\newcommand{\qed}{\nobreak \ifvmode \relax \else
      \ifdim\lastskip<1.5em \hskip-\lastskip
      \hskip1.5em plus0em minus0.5em \fi \nobreak
      \vrule height0.5em width0.5em depth0.25em\fi}

\begin{document}
\bibliographystyle{amsplain}
\title{ A note on the value distribution of \\
$f^l(f^{(k)})^n$}

\author{Yan Jiang and Bin Huang
}

%\author{Bin Huang%
 % \thanks{E-mail address: \texttt{huangbincscu@163.com \\
%Mathematics Subject Classifications (2010): Primary 30D35}}}
%  \affil{}

\date{ }

\maketitle \baselineskip 13pt

\begin{minipage}{136mm}
\abstract{Let $f$ be a transcendental meromorphic function in the
complex plane $\mathbb{C}$, and $a$ be a nonzero complex number . We
give quantitative estimates for the characteristic function $T(r,f)$
in terms of $N(r,1/(f^l(f^{(k)})^n-a))$, for integers $k$, $l$, $n$
greater than 1.
We conclude that $f^l(f^{(k)})^n$ assumes every nonzero finite value infinitely often.}
\end{minipage}
\section{Introduction}

Let $f$ be a transcendental meromorphic function in the complex
plane $\mathbb{C}$. In this article, we use the standard notations
in the sense of Nevanlinna \cite{W.K.Hayman}, such as $T(r, f)$,
$N(r, f)$, $\bar{N}(r, f)$, $m(r, f)$, $S(r, f)$, $\delta(a, f)$. In
particular, $T(r,f)$ the characteristic function and $\bar{N}(r,f)$
is a counting function with respect to poles of $f$, ignoring
multiplicities. We shall use the symbol $S(r,f)$ to denote an error
term $v(r)$ satisfying $T(r,v(r))=o\left(T(r,f)\right)$ as $r\rightarrow
\infty$, possibly outside a set of finite linear measure. Throughout
this paper a small function (with respect to $f$) means a function
$\varphi(z)$ meromorphic in $\mathbb{C}$ satisfying
$T(r,\varphi)=S(r,f)$. In addition, in this paper, we use another
type of small function $S^*(r,f)$ which has the property
$S^*(r,f)=o\left(T(r,f)\right)$ as $r\rightarrow \infty$, $r\not\in
E$, where $E$ is a set of logarithmic density 0.

A meromorphic function $f$ is rational if and only if $T(r,f)=O(\log
r)$ (see \cite{Goldberg}). The quantity
$$
\delta(a,f)=\mathop {\lim \inf}\limits_{r \to \infty }
\frac{m(r,1/(f-a))}{T(r,f)}=1-\mathop {\lim \sup}\limits_{r \to
\infty } \frac{N(r,1/{(f-a)})}{T(r,f)}
$$
is called the deficiency of $f$ at the point $a$. Another deficiency is defined by
$$
\Theta(a,f)=1-\mathop {\lim \sup}\limits_{r \to \infty }
\frac{\bar{N}(r,1/{(f-a)})}{T(r,f)}.
$$
Note that $0\leq\delta(a,f)\leq\Theta(a,f)\leq 1$.

The First Fundamental Theorem of the value distribution theory
 due to Nevanlinna is utilized frequently in this note. It is stated as the following property:
$$
T\left(r,\frac{1}{f-a}\right)=T(r,f)+O(1)
$$
for any constant $a$. The details can be found in \cite{Goldberg}
for example. A root of the equation $f(z)=a$ $(1/f(z)=0$ for
$a=\infty)$ will be called an $a-$point of the function $f(z)$ for
$a\in \mathbb{C}\cup \{\infty\}$. $a$ is called a Picard exceptional
 value of a function $f(z)$ if the number of its $a-$points in $\mathbb{C}$ is finite.

The aim of this paper is to look for estimates with respect to
$f^l(f^{(k)})^n$. The following well-known estimate is due to Hayman
\cite[Theorem 3.5]{W.K.Hayman}.
\begingroup
\setcounter{tmp}{\value{theo}}% store current value of theorem counter
\setcounter{theo}{0} %assign desired value to theorem counter
\renewcommand\thetheo{\Alph{theo}}% locally redefine the representation of the theorem counter
\begin{theo}
Let $f$ be a meromorphic and
transcendental function in the plane, $l$ be a positive integer, and
$a$, $b$ be constants with $b\neq 0$. Then
\begin{equation}\label{eq00}
T(r,f)\leq
\left(2+\frac{1}{l}\right)N\left(r,\frac{1}{f-a}\right)+\left
(2+\frac{2}{l}\right)\bar{N}\left(r,\frac{1}{f^{(l)}-b}\right)+S(r,f).
\end{equation}
\end{theo}

%\begingroup
%\setcounter{ccc}{\value{corollary}}% store current value of theorem counter
%\setcounter{corollary}{0} %assign desired value to theorem counter
%\renewcommand\thecorollary{\Alph{corollary}}% locally redefine the representation of the theorem counter

Hayman also concluded a corollary from the previous inequality.
\begin{corollary-non} Under the same
assumptions as in Theorem A, either $f$ assumes every finite value
infinitely often or $f^{(l)}$ assumes every finite value except
possibly zero infinitely often.
\end{corollary-non}

%\endgroup
Moreover, Hayman conjectured that if $f$ is a transcendental
meromorphic function and $l\geq 1$, then $f^lf'$ takes every finite
nonzero value infinitely often.
%In 1969, Sons (\cite{Sons}, Theorem 2) proved similar results for
%the case $n\geq 2$, however under the additional assumptions
%$n_k\geq 1$ and
%$$
%2^k\left(2n_0+\sum_{i=0}^k(1+i)n_i\right)<(2^k+2n_0-1)\left(\sum_{i=0}^k(1+i)n_i)
%\right).
%$$
This conjecture has been confirmed by himself in \cite{W.K.Hayman}
for $l\geq 3$, by Mues \cite{Mues} for $l=2$ and by Bergweiler and
Eremenko \cite{Ber1} for $l=1$. During the past decades, a sequence
of related research have been made. In 1982, Doeringer
\cite[Corollary 1]{Doeringer} proved that for a transcendental
meromorphic function $f$, the only possible
 Picard exceptional value is zero for a differential monomial $f^l(f^{(k)})^n$
 when $l\geq 3$.
 In 1994, Tse and Yang \cite{C.K.Tse} gave an estimate of $T(r,f)$ for $l=1$ and $l=2$ and confirmed  the only possible
 Picard exceptional value is zero. In 1996, Yang
and Hu \cite[Theorem 2]{Yang} proved that if
$\delta(0,f)>3/(3(l+n)+1)$ with positive integers $k$, $l$, $n$,
then for a nonzero finite complex number
 $a$, $f^l(f^{(k)})^n-a$ has infinitely many zeros. In 2002, Li and
Wu \cite{Li2} obtained that for a nonzero finite complex number
 $a$ and  positive integers $l$, $k$ with $l\geq 2$, there exists a constant $M > 0$ such that
$$T(r, f) < M \bar{N}\left(r,\frac{1}{f^lf^{(k)}-a}\right)+ S(r, f).$$
 In 2003, Wang \cite{Wang} studied the zeros of $f^lf^{(k)}-\phi$ for a small
meromorphic function $\phi(z)\not\equiv 0$, and
 verified that for $l\geq 2$, $f^lf^{(k)}-\phi$ had infinitely many zeros if the poles of $f$
were multiple. In 2004, Alotaibi \cite{Alotaibi} gave an estimate
and showed that the function $f(f^{(k)})^n-\phi$ has infinitely many
zeros for a small function $\phi(z)\not\equiv 0$, when $n\geq 2$.

 We introduce a result given by Lahiri and Dewan \cite[Theorem 3.2]{Lahiri1}.

\begin{theo}\label{lahiri}
Let $f$ be a transcendental meromorphic function and $a$, $\alpha$
be both small functions of $f$ without being identically to zero and
infinity. If $\psi=\alpha f^l(f^{(k)})^n$, where $l(\geq 0)$,
$n(\geq 1)$, $k(\geq 1)$ are integers, then
\begin{eqnarray}\label{eqlahiri}
(l+n)T(r,f)\leq \bar{N}\left(r,f\right)+
\bar{N}\left(r,\frac{1}{f}\right)
+nN_{(k)}\left(r,\frac{1}{f}\right)+\bar{N}\left(r,\frac{1}{\psi-a}\right)
+S(r,f),
\end{eqnarray}
where $N_{(k)}(r,1/f)$ is the counting function of zeros of $f$ with
multiplicity $q$ counted $\min \{q,k\}$ times.
\end{theo}

\noindent {\it Remark.} Inequality (\ref{eqlahiri}) implies that for
$l\geq 3$, $n\geq 1$, $k\geq 1$,
\begin{equation}\label{eq0137}
T(r,f) \leq \frac{1}{l-2}N\left(r,\frac{1}{f^l(f^{(k)})^n -
a}\right) + S(r,f),
\end{equation}
then
\begin{eqnarray}\label{eqla}
\delta(a,\psi)\leq \Theta(a,\psi)\leq 1-\frac{l-2}{nk+n+l}.
\end{eqnarray}

 However, this result is still worth refining. In the current paper, we obtained an estimate
corresponding to the case $k$, $l$, $n$ all greater than 1, and
in our proof, we use a very important inequality of Yamanoi.

\begin{theo}\textbf{\emph{\cite[Yamanoi]{Yamanoi}}}\label{th000}
Let $f$ be a meromorphic and transcendental function in the complex
plane and let $k\geq 2$ be an integer, $A\subset \mathbb{C}$ be a finite set of complex numbers. Then we
have
\begin{equation}\label{eq000}
(k-1)\bar{N}(r,f)+\sum_{a\in A}N_1\left(r,\frac{1}{f-a}\right)\leq
N\left(r,\frac{1}{f^{(k)}}\right)+S^*(r,f),
\end{equation}
where
$$
N_1\left(r,\frac{1}{f-a}\right)=N\left(r,\frac{1}{f-a}\right)-
\bar{N}\left(r,\frac{1}{f-a}\right).
$$
\end{theo}
\endgroup

\noindent {\it Remark.} Actually, this theorem is obviously a correspondent result to the
famous Gol'dberg Conjecture, which says that for a transcendental
meromorphic function $f$ and $k\geq 2$, then $
\bar{N}(r,f)\leq N\left(r,1/{f^{(k)}}\right)+S(r,f)$. There
is obviously a special case of Yamanoi's inequality when $A$ is an
empty set. The following special case is the one we use in our proof,
\begin{equation}\label{eq001}
(k-1)\bar{N}(r,f)\leq N\left(r,\frac{1}{f^{(k)}}\right)+S^*(r,f).
\end{equation}

 In this paper, we continue to consider the
general form $f^l(f^{(k)})^n-a$  for a nonzero constant $a$.
The following theorem improved Theorem \ref{lahiri} in some sense.
\setcounter{theo}{\thetmp} % restore value of theorem counter
\begin{theo}\label{th1} Let $f$ be a transcendental meromorphic
function in $\mathbb{C}$, $l$, $n$, $k$ be integers greater than 1 and $a$
be a nonzero constant. Then
\begin{equation}\label{eq0138}
T(r,f)\leq
\frac{1}{l-1}N\left(r,\frac{1}{f^l(f^{(k)})^n-a}\right)+S^*(r,f).
\end{equation}
\end{theo}
\noindent {\it Remark.} If the differential monomials
$f^l(f^{(k)})^n$ is allowed to take $l\geq 2$, $n\geq 2$, $k\geq 2$, then (\ref{eq0138}) is better than (\ref{eq0137}) with the except of a finite set of logarithmic density 0. If
the case $k=1$, $l\geq 3$ or $n=1$, $l\geq 3$ occurs,
(\ref{eq0137}) might be the best choice so far.

Another important remark should be made here. As we realized that
for general form $f^l(f^{(k)})^n$, except the cases stated in
Theorem \ref{lahiri} and Theorem \ref{th1}, two cases are inevitably
excluded: $l=1$, $n\geq 1$, $k\geq 1$ and $l=2$, $n\geq 1$, $k\geq
1$. We summarize the known estimates of these two cases. For the
case $l=2$, $n=k=1$, Zhang \cite{Zhang} obtained a quantitative
result, proving that the inequality
$T(r,f)<6N(r,1/(f^2f'-1))+S(r,f)$ holds. For the case $l=2$, $n=1$,
$k>1$ the the inequality is due to Huang and Gu \cite{Huang}. For
the case $l=1$, $n\geq 2$, $k\geq 1$,
 by Li and Yang \cite{Li} and
Alotaibi \cite{Alotaibi} gave two different inequalities for the estimates independently.
For the case $l=n=1$, $k\geq 1$, again
Alotaibi \cite{Alotaibi2} obtained an estimate provided that $
N_{1)}\left(r,\frac{1}{f^{(k)}}\right)=S(r,f)$, where
$N_{1)}\left(r,\frac{1}{f^{(k)}}\right)$ is the counting function of
simple zeros of $f^{(k)}$, as well, Wang \cite{Wang} gave an
estimate but under the additional condition that multiplicities of all
poles of $f$ are at least $3$ and $N_{1)}(r,1/f)\leq \lambda
T(r,f)$, where $\lambda<1/3$ is a constant.

Though these cases are excluded in Theorem \ref{th1}, our estimate
is considered to be stronger compared to the known results so far. Furthermore, it is natural to estimate the
deficiency of $f^l(f^{(k)})^n$ by making use of Theorem \ref{th1}. This leads us to the following.

\begin{theo}\label{th2}
Let $f$ be a transcendental meromorphic function in $\mathbb{C}$,
$k$, $l$, $n$ be positive integers all greater than 1 and $a$ be a nonzero constant.
Then
$$
\delta(a,f^l(f^{(k)})^n)\leq 1-\frac{l-1}{nk+n+l}.
$$
%When $k=1$, $l\geq 3$, $n\geq 1$,
%$$
%\delta(a,f^l(f')^n)\leq 1-\frac{l-2}{2n+l}.
%$$
\end{theo}
\noindent {\it Remark.} Since for a nonzero constant $a$,
$\delta(a,f^l(f^{(k)})^n)< 1$, Theorem \ref{th2} also implies that
the possible Picard exceptional value of $f^l(f^{(k)})^n$
 is zero for $k\geq 2$, $l\geq 2$, $n\geq 2$. We like to state these results as a corollary here.

\begin{corollary}
Under the same conditions as Theorem \ref{th1}, $f^l(f^{(k)})^n$
assumes every finite value except possibly zero infinitely often.
\end{corollary}
\noindent {\it Remark.}  In fact, this kind of result is not brand
new. There are already a couple of known results implying that for
any positive integers $k$, $l$, $n$, the function $f^l(f^{(k)})^n$
assumes every finite value except possibly zero infinitely often.
The readers should see Lahiri and Dewan \cite{{Lahiri1},{Lahiri2}},
Steinmetz \cite{Ste}, Wang \cite{Wang1}, Alotaibi
\cite{{Alotaibi},{Alotaibi2}} and Li and Wu \cite{Li2} for further
details.

%Bergweiler and Eremenko \cite{Ber1}  show that $f'f$ assumes every
%finite nonzero complex value infinitely often for a transcendental
%meromorphic function $f$. W. Bergweiler \cite{Ber2} also concluded
%that if transcendental meromorphic function $f$ is of finite order
%and $c\not\equiv 0$ is a polynomial, then $f'f-c$ has infinitely
%many zeros.
\bigskip

Lemmas used for the proof of Theorem 1.1 and Theorem
1.2 are presented in Section 2. The proofs of Theorem \ref{th1} and
Theorem \ref{th2} are placed in Section 3 and 4 respectively. In the
last section, we give an application to the sum of deficiencies.

\section{Lemmas}

\quad Before we proceed to the proofs of the theorems, we need the
following lemmas.

\begin{lemma}\textbf{\emph{\cite[Theorem 3.1]{W.K.Hayman}}}\label{lm1} \it{Let $f$ be a non-constant
meromorphic function in the complex plane, $l$ be positive
integer, $a_0(z)$, $a_1(z)$,$\cdots$, $a_l(z)$ be meromorphic
functions in the plane satisfying $T\left(r,a_\nu(z)\right)=S(r,f)$ for $\nu=0$, $1$, $\ldots$, $l$ (
as $r\rightarrow +\infty$) and
$$
\psi(z)=\sum_{\nu=0}^l a_\nu(z)f^{(\nu)}(z).
$$
Then
$$
m\left(r,\frac{\psi}{f}\right)=S(r,f).
$$}
\end{lemma}

In particular, this lemma implies $m\left(r, f^{(l)}/f\right)= S(r,f)$ and
$m\left(r,f^{(l+1)}/f^{(l)}\right)= S(r,f^{(l)})$.
\bigskip

\begin{lemma}\textbf{\emph{\cite[p. 99]{Goldberg}}}\label{lm2} \it{ Let $f$ be a non-constant
meromorphic function in the complex plane, $k$ be a positive
integer.
Then
\begin{eqnarray}\label{eq037}
T(r,f^{(k)})\leq (k+1)T(r,f)+S(r,f).
\end{eqnarray}}
\end{lemma}

In particular, $S(r,f^{(k)})\leq S(r,f)$. Inequality (\ref{eq037}) will be used often in this note without reference.

\begin{lemma}\label{lm3}
Let $f$ be a transcendental meromorphic function in the plane. Then the differential monomial
$$
\psi= f^{l}(f^{(k)})^n
$$
is transcendental, where $l$, $n$ and $k$ are positive integers.
\end{lemma}

\noindent {\it Proof.}
Since we have
$$
\frac{1}{f^{l+n}}=\left(\frac{f^{(k)}}{f}\right)^n\frac{1}{\psi}.
$$
We obtain from Lemma \ref{lm1} and the First Fundamental Theorem
that
\begin{eqnarray}\label{eq1111}
(l+n)T(r,f) &\leq & nT\left(r,\frac{f^{(k)}}{f}\right)+T\left(r,\frac{1}{\psi}\right)\nonumber\\
&\leq& n N\left(r,\frac{f^{(k)}}{f}\right)+T\left(r,\frac{1}{\psi}\right)+S(r,f)\nonumber\\
&\leq& nk\left[\bar{N}(r, f)+\bar{N}\left(r, \frac{1}{f}\right)\right]+T\left(r,\frac{1}{\psi}\right)+S(r,f).
 \end{eqnarray}

  Since $\bar{N}(r, f)\leq \bar{N}(r, \psi)+S(r,f)$ and $\bar{N}\left(r, \frac{1}{f}\right) \leq \bar{N}\left(r, \frac{1}{\psi}\right)+S(r,f)$,
we can simplify inequality (\ref{eq1111}) to
$$
(l+n)T(r,f)  \leq (2nk+1)T\left(r,\frac{1}{\psi}\right)+S(r,f).
$$
Because $f$ is transcendental, we conclude that $\psi$ is
transcendental.\qed

%\begin{lemma}\label{lm4} Let $f$ be a transcendental meromorphic function in $\mathbb{C}$ and $k$, $l$, $n$
%be positive integers. Then $f^l(f^{(k)})^n$ is non-constant.
%\end{lemma}

%\noindent {\it Proof.} Suppose that $f^l(f^{(k)})^n $
%is a constant $C$. It is clear that $C \ne 0$, so $f$ has no zeros. Thus, $
%C /f^{l + n} = \left((f^{(k)}) / f\right)^n$. By Lemma 2.1 and the First Fundamental Theorem,
%\begin{eqnarray*}
% T(r,f^{l + n}) &=& T\left(r,\frac{C}{f^{l + n}}\right) + O(1) \\
% &= & m\left(r,\frac{C}{f^{l+n}}\right)+ O(1) \\
 %&= & m\left(r,\left(\frac{f^{(k)}}{f}\right)^{n}\right)+ O(1) \\
 %&= & S(r,f).
 %\end{eqnarray*}
%On the other hand
%$$T(r,f^{l+n})=(l+n)T(r,f)+O(1),
%$$
%which implies that $T(r,f)=S(r,f)$. However, this is impossible by assumption.\qed
\bigskip
\begin{lemma}\label{lm11} Let $f$ be a transcendental meromorphic function in $\mathbb{C}$,
let $k$, $l$, $n$ be positive integers, and set
$$g =f^l(f^{(k)})^n - 1.$$
Then,
$$T(r,g)\leq O\left(T(r,f)\right),
$$
as $r\rightarrow \infty$, possibly outside a set of finite linear
measure.
\end{lemma}

\noindent{\it Proof.} Note that $N\left(r,f^l(f^{(k)})^n\right)=O(N(r,f))$ and
 $m\left(r,f^{(k)}/f\right)=S(r,f)$ by Lemma \ref{lm1}. Applying the First Fundamental Theorem, we get
\begin{eqnarray*}
T(r,g) &=& T\left(r,f^l(f^{(k)})^n-1\right)\\
&=& N\left(r,f^l(f^{(k)})^n\right)+m\left(r,f^l(f^{(k)})^n\right)+O(1)\\
&\leq& O\left(N(r,f)\right)+lm(r,f)+nm\left(r,f^{(k)}\right)++O(1)\\
&\leq& O\left(N(r,f)\right)+lm(r,f)+nm(r,f)+nm\left(r,\frac{f^{(k)}}{f}\right)+O(1)\\
&=& O\left(T(r,f)\right)+S(r,f).
\end{eqnarray*}

We can see that
$$T(r,g')\leq N(r,g')+m(r,g)+S(r,g)\leq
T(r,g)+S(r,g).$$
 Hence
\[ T(r,g)\leq O\left(T(r,f)\right).\qed\]

\section{Proof of Theorem \ref{th1}.}
Without loss of generality, we assume $a=1$. $g =f^l(f^{(k)})^n -
1$.
 By Lemma \ref{lm3}, we
know that $g$ is not constant. Since
$$\frac{1}{f^{l + n}} = \left(\frac{f^{(k)}}{f}\right)^n -
\frac{g'}{f^{l + n}}\left(\frac{g}{g'}\right),
$$
it follows that
$$
m\left(r,\frac{1}{f^{l+n}}\right)\leq
m\left(r,\frac{g}{g'}\right)+m\left(r,\frac{g'}{f^{l+n}}\right)+S(r,f).
$$

Note that
$$\frac{{g}'}{f^{l + n}} =
l\frac{{f}'}{f}\left(\frac{f^{(k)}}{f}\right)^n + n\frac{f^{(k +
1)}}{f}\left(\frac{f^{(k)}}{f}\right)^{n - 1},
$$
which implies
$$
m\left(r,\frac{g'}{f^{l+n}}\right)=S(r,f).
$$
 Therefore, we have
$$
m\left(r,\frac{1}{f^{l+n}}\right)\leq
m\left(r,\frac{g}{g'}\right)+S(r,f).
$$

We know that the poles of $g'/g$ come from the zeros and poles of
$g$, and all are simple. The poles of $g/g'$ come from zeros of $g'$
which are not zeros of $g$, preserving multiplicity. Hence, we get
\begin{equation}
\label{eq27}
N\left(r,\frac{g'}{g}\right)=\bar{N}\left(r,\frac{1}{g}\right)+\bar{N}(r,g),
\end{equation}
and
\begin{equation}
\label{eq28}
N\left(r,\frac{g}{g'}\right)=N\left(r,\frac{1}{g'}\right)-\left(N\left(r,\frac{1}{g}\right)-\bar{N}\left(r,\frac{1}{g}\right)\right).
\end{equation}
By combining (\ref{eq27}) with (\ref{eq28}),
\begin{equation}
\label{eq29}
N\left(r,\frac{g'}{g}\right)-N\left(r,\frac{g}{g'}\right)=\bar{N}(r,g)+N\left(r,\frac{1}{g}\right)-N\left(r,\frac{1}{g'}\right).
\end{equation}

By Lemma \ref{lm11}, we know that
$$m(r,g'/g)=S(r,g)\leq S(r,f), \quad\quad
\bar{N}(r,g)=\bar{N}(r,f).$$
Applying the First Fundamental Theorem
and (\ref{eq29}),
\begin{eqnarray*}
m\left(r,\frac{1}{f^{l+n}}\right)
&=& (l + n)m\left(r,\frac{1}{f}\right)\\
&\leq& N\left(r,\frac{g'}{g}\right)- N\left(r,\frac{g}{g'}\right) + m\left(r,\frac{g'}{g}\right)+ S(r,f) \\
&\leq& N\left(r,\frac{g'}{g}\right) - N\left(r,\frac{g}{g'}\right) +
S(r,g)+ S(r,f)
\end{eqnarray*}
\begin{equation}
\label{eq5} \ \ \ \ \ \ \ \ \ \ \ \ \ \ =\
\bar{N}(r,f)+N\left(r,\frac{1}{g}\right) -
N\left(r,\frac{1}{g'}\right) +S(r,f).
\end{equation}
Here we add $N(r,1/f^{l+n})$ to both sides of inequality
(\ref{eq5}), then
\begin{equation}\label{eq16}
(l + n)T\left(r,\frac{1}{f}\right)\leq\bar {N}(r,f) +
N\left(r,\frac{1}{g}\right) - N\left(r,\frac{1}{g'}\right) +
N\left(r,\frac{1}{f^{l+n}}\right)+ S(r,f).
\end{equation}
Note that
$g'=f^{l-1}\left(f^{(k)}\right)^{n-1}\left(lf^{(k)}f'+nff^{(k+1)}\right)$,
which implies
\begin{equation}
\label{eq17}
(l-1)N\left(r,\frac{1}{f}\right)+(n-1)N\left(r,\frac{1}{f^{(k)}}\right)\leq
N\left(r,\frac{1}{g'}\right).
\end{equation}
Substituting (\ref{eq17}) into (\ref{eq16}), we get
\begin{eqnarray*}
T\left(r,\frac{1}{f^{l+n}}\right) &\leq&
\bar{N}(r,f)+N\left(r,\frac{1}{g}\right)+N\left(r,\frac{1}{f^{l+n}}\right)-(l-1)N\left(r,\frac{1}{f}\right)\\
&&-(n-1)N\left(r,\frac{1}{f^{(k)}}\right)+S(r,f).
\end{eqnarray*}
Hence,
\begin{equation}
\label{eq172} (l+n)T(r,f)\leq
N\left(r,\frac{1}{g}\right)+\bar{N}(r,f)+(n+1)N\left(r,\frac{1}{f}\right)-(n-1)N\left(r,\frac{1}{f^{(k)}}\right)+S(r,f).
\end{equation}

Inequality (\ref{eq001}) implies that for $k\geq 2$,
\begin{eqnarray}
\label{eq173} (k-1)\bar{N}(r,f) \leq
N\left(r,\frac{1}{f^{(k)}}\right)+S^*(r,f).
\end{eqnarray}
Now by combining inequality (\ref{eq172}) and (\ref{eq173}), we have
\begin{eqnarray*}
(l+n)T(r,f)&\leq&
N\left(r,\frac{1}{g}\right)+\frac{1}{k-1}N\left(r,\frac{1}{f^{(k)}}\right)+(n+1)N\left(r,
\frac{1}{f}\right)\\
 &&-(n-1)N\left(r,\frac{1}{f^{(k)}}\right)+S(r,f)+S^*(r,f).
\end{eqnarray*}

Since $1/{(k-1)}-n+1\leq 0$ for $n\geq 2$, $k\geq 2$, then
\begin{eqnarray}\label{eq123}
(l+n)T(r,f)\leq
N\left(r,\frac{1}{g}\right)+\left(n+1\right)N\left(r,\frac{1}{f}\right)+S^*(r,f).
\end{eqnarray}
Since $l-1>0$ and $\left(n+1\right)N\left(r,1/f\right)\leq
\left(n+1\right)T\left(r,f\right)$, then
\begin{eqnarray}\label{eq008}
T(r,f)\leq
\frac{1}{l-1}N\left(r,\frac{1}{f^l(f^{(k)})^n-1}\right)+S^*(r,f).
\end{eqnarray}
Replacing the number 1 in $f^l(f^{(k)})^n-1$ by any nonzero constant
$a$, the inequality (\ref{eq0138}) is obtained. The proof is
completed.\qed

% On the other hand, note that
%\begin{equation}
%\label{eq8} (l + n)T(r,f) = (l - 1)(N(r,f) + m(r,f)) + (n +
%1)\left(N\left(r,\frac{1}{f}\right) + m\left(r,\frac{1}{f}\right)\right) + O(1).
%\end{equation}
%Combining the inequalities (\ref{eq7}) and (\ref{eq8}),
%\begin{eqnarray*}(l-1)N(r,f)-\bar{N}(r,f)+N\left(r,\frac{1}{h}\right)-(n+1)m\left(r,\frac{1}{f}\right)
%&\leq& N\left(r,\frac{1}{g}\right)-(l-1)m(r,f)+S(r,f)\\
%&\leq& N\left(r,\frac{1}{g}\right)+S(r,f).
%\end{eqnarray*}
%Thus, the inequality (\ref{eq4}) is concluded.

\section{Proof of Theorem \ref{th2}.}

Set $\psi=f^l(f^{(k)})^n$. Inequality (\ref{eq008}) is stated that
\begin{equation}
\label{eq072} T(r,f)\leq
\frac{1}{l-1}N\left(r,\frac{1}{\psi-a}\right)+S^*(r,f)
\end{equation}
for $l\geq 2$, $n\geq 2$, $k\geq 2$.
By the definition of
$\delta(a,f)$ and the First Fundamental Theorem, we obtain
\begin{eqnarray}\label{eq111}
T(r,\psi)&\leq& (nk+n+l)T(r,f)+S(r,f)\nonumber\\
&\leq& \frac{nk+n+l}{l-1} N\left(r,\frac{1}{\psi-a}\right)+S^*(r,f).
\end{eqnarray}
By combining two inequalities (\ref{eq072}) and (\ref{eq111}), we have
$$
N\left(r,\frac{1}{\psi-a}\right)\geq
\frac{l-1}{nk+n+l}T(r,\psi)-S^*(r,f).
$$

Since
\begin{eqnarray*}
T(r,f)&=&\frac{1}{l}T(r,f^l)\leq T\left(r,(f^{(k)})^n\right)+T(r,\psi)\\
&\leq& O\left(T\left(r,\psi\right)\right),
\end{eqnarray*}
then we deduce that
$$
\mathop {\lim \inf}\limits_{r \to \infty}
\frac{S^*(r,f)}{T(r,\psi)}=\mathop {\lim \inf}\limits_{r \to \infty}
\frac{S^*(r,f)}{T(r,f)}\frac{T(r,f)}{T(r,\psi)}= 0.
$$
Therefore, by the definition of deficiency,
\begin{eqnarray*}
\delta(a,\psi) &=& 1-\mathop {\lim \sup}\limits_{r \to
\infty}\frac{N\left(r,\frac{1}{\psi-a}\right)}{T(r,\psi)}\\
 &\leq& 1-\mathop {\lim
\sup}\limits_{r \to
\infty}\frac{\frac{l-1}{nk+n+l}T(r,\psi)-S^*(r,f)}{T(r,\psi)}\\
&\leq& 1-\frac{l-1}{nk+n+l}+\mathop {\lim
\inf}\limits_{r \to
\infty}\frac{S^*(r,f)}{T(r,\psi)}\\
&=& 1-\frac{l-1}{nk+n+l}.\qed
\end{eqnarray*}

%%In the sense of the definition of $\delta(a,f)$ and by the First
%Fundamental Theorem, we obtain
%\begin{eqnarray*}
%T(r,\psi)&\leq& (l+n+nk)T(r,f)+S(r,f)\\
%&\leq& \frac{l+n+nk}{l-1-\frac{1+\varepsilon}{k}} N\left(r,\frac{1}{\psi-a}\right)+S(r,f).
%\end{eqnarray*}
%$$N\left(r,\frac{1}{\psi-a}\right)\geq \frac{l-1-\frac{1+\varepsilon}{k}}{l+n+nk}T(r,\psi)-S(r,f).$$
%We know that $S(r,f)\leq S(r,\psi)$ by Theorem 1.3. On the other hand, by Lemma 2.3, $S(r,f)\geq S(r,g)=S(r,\psi)$. We deduce that
%$$\mathop {\lim \inf}\limits_{r \to \infty} \frac{S(r,f)}{T(r,\psi)}= 0,$$
%possibly except a set of finite linear measure. Hence, we have
%\begin{eqnarray*}
%\delta(a,\psi) &=& 1-\mathop {\lim \sup}\limits_{r \to
%\infty}\frac{N\left(r,\frac{1}{\psi-a}\right)}{T(r,\psi)}\\
% &\leq& 1-\mathop {\lim
%\sup}\limits_{r \to
%\infty}\frac{\frac{l-1-\frac{1+\varepsilon}{k}}{l+n+nk}T(r,\psi)-S(r,f)}{T(r,\psi)}\\
%&\leq& 1-\frac{l-1-\frac{1+\varepsilon}{k}}{l+n+nk}+\mathop {\lim \inf}\limits_{r \to
%\infty}\frac{S(r,\psi)}{T(r,\psi)}\\
%&=& 1-\frac{l-1-\frac{1+\varepsilon}{k}}{l+n+nk}.
%\end{eqnarray*}
\section{An application}
After Yamanoi's result was published in 2013, there are some results
about deficieny relations came out by using his important theorem.
We take a result from Fang and Wang \cite{Fang1} as a good example
here, and we analogue their steps to get an estimate of the sum of
deficiencies of $f^l(f^{(k)})^n$.

\begingroup
\setcounter{tmp}{\value{theo}}% store current value of theorem counter
\setcounter{theo}{3} %assign desired value to theorem counter
\renewcommand\thetheo{\Alph{theo}}% locally redefine the representation of the theorem counter
\begin{theo}\textbf{\emph{\cite[Propostion 2]{Fang1}}}
\it{Let $f$ be a meromorphic and
transcendental function in the complex plane, $k$ be a positive integer, and
$P$ be the set of all polynomials. Then
\begin{equation*}
\sum_{b\in \mathbb{C}} \delta{(b,f^{(k)})}\leq 1-(k-1)\left(1-\Theta_E\left(\infty,f^{(k)}\right)\right),
\end{equation*}
where for $r\not\in E$,
$$
\Theta_E\left(\infty,f^{(k)}\right)=1-\mathop
{\lim\sup}\limits_{r\to \infty}
\frac{\bar{N}(r,f^{(k)})}{T(r,f^{(k)})},
$$
where $E(=M(K)\cup E_1)$, $M(K)$ is a set of finite upper
logarithmic density and $E_1$ is a set of finite measure.}
\end{theo}
\endgroup
We need the following lemma for our calculation. This lemma is as
well used in paper \cite{Fang1}.
\begin{lemma}\textbf{\emph{\cite[p. 33]{W.K.Hayman}}}\label{lmp} \it{Let $a_1$, $a_2$, $\cdots$, $a_q$, where $q>2$, be distinct finite
 complex numbers. Then
 $$\sum_{i=1}^q m\left(r, \frac{1}{f-a_i}\right)
 \leq m\left(r,\sum_{i=1}^{q}\frac{1}{f-a_i}\right)+O(1).$$}
\end{lemma}

\bigskip
\begin{theo}\label{th3}
Let $f$ be a transcendental meromorphic function in $\mathbb{C}$,
$k$, $l$, $n$ be positive integers all at least $2$ and $a_i\in \mathbb{C}$ be constants, $i=1,2,\cdots,q$.
Then
$$
\sum_{i=1}^q\delta(a_i,f^l(f^{(k)})^n)\leq  1+\frac{1}{nk+n+l}.
$$
\end{theo}

\noindent {\it Proof}. By Nevanlinna theory, for constants $a_i\in \mathbb{C}$, the sum of deficiencies of function $f$  are defined by
\begin{eqnarray}\label{eqd}
\sum_{i=1}^{q} \delta{(a_i,f)}=\mathop {\lim \inf}\limits_{r \to
\infty}\sum_{i=1}^{q}\frac{m\left(r, \frac{1}{f-a_i}\right)}{T\left(r,f\right)}.
\end{eqnarray}
Let $\psi=f^l(f^{(k)})^n$. By Lemma \ref{lmp}, we have
\begin{eqnarray}\label{eqs}
\sum_{i=1}^q m\left(r, \frac{1}{\psi-a_i}\right)
&\leq& m\left(r,\sum_{i=1}^{q}\frac{1}{\psi-a_i}\right)+O(1)\nonumber\\
&\leq& m\left(r,\sum_{i=1}^{q}\frac{\psi''}{\psi-a_i}\right)+m\left(r,\frac{1}{\psi''}\right)+S(r,f)\nonumber\\
&\leq& T(r,\psi'')-N\left(r,\frac{1}{\psi''}\right)+S(r,f)\nonumber\\
&\leq& N(r,\psi'')+m(r,\psi'')-N\left(r,\frac{1}{\psi''}\right)+S(r,f)\nonumber\\
&\leq& N(r,\psi)+2\bar{N}(r,\psi)+m(r,\psi)-N\left(r,\frac{1}{\psi''}\right)+S(r,f).
\end{eqnarray}
By Yamanoi's result (\ref{eq001}), it follows from inequality (\ref{eqs}) that
\begin{eqnarray}\label{eqp}
\sum_{i=1}^q m\left(r, \frac{1}{\psi-a_i}\right)
&\leq& T(r,\psi)+2\bar{N}(r,\psi)-\bar{N}\left(r,\psi\right)+S(r,f)\nonumber\\
&\leq& T(r,\psi)+\bar{N}(r,f)+S(r,f)\nonumber\\
&\leq& T(r,\psi)+T(r,f)+S(r,f).
\end{eqnarray}
By Theorem  \ref{th1} and Theorem  \ref{th2}, it follows from inequality (\ref{eqp}) that,
\begin{eqnarray*}
\sum_{i=1}^q m\left(r, \frac{1}{\psi-a_i}\right)&=&\mathop {\lim \inf}\limits_{r \to
\infty}\sum_{i=1}^{q}\frac{m\left(r, \frac{1}{\psi-a_i}\right)}{T\left(r,\psi\right)}\\
&\leq& 1+\mathop {\lim \inf}\limits_{r \to
\infty}\frac{T(r,f)}{T(r,\psi)}+S(r,f)\\
&\leq& 1-\frac{1}{l-1}\left(1-\mathop {\lim \sup}\limits_{r \to
\infty}\frac{N\left(r,\frac{1}{\psi-a}\right)}{T(r,\psi)}-1\right)\\
&=& 1-\frac{1}{l-1}\left(\delta\left(a,\psi\right)-1\right)\\
&\leq& 1-\frac{1}{l-1}\left( 1-\frac{l-1}{nk+n+l}-1\right)\\
&=& 1+\frac{1}{nk+n+l}.\qed
\end{eqnarray*}

\section{Acknowledgement}
\quad The authors are very grateful to Professor Toshiyuki Sugawa and every member in the seminars for their valuable
suggestions and comments, which helped a lot in improving the paper.

%\bibliography{jiang2013}
\bigskip
\providecommand{\bysame}{\leavevmode\hbox
to3em{\hrulefill}\thinspace}
\providecommand{\MR}{\relax\ifhmode\unskip\space\fi MR }
\providecommand{\MRhref}[2]{%
  \href{http://www.ams.org/mathscinet-getitem?mr=#1}{#2}
} \providecommand{\href}[2]{#2}

\vskip 1cm
\noindent
\large{Yan Jiang}\\
Graduate School of Information
Sciences,\\ Tohoku University, Sendai 980-88579, JAPAN\\
\bigskip
E-mail address: jiang@ims.is.tohoku.ac.jp; gesimy038@tom.com\\
\large{Bin Huang}\\
Department of Mathematics and Computing Science,\\
Changsha University of Science and Technology, Changsha 410076,
P.R.China\\
E-mail address: huangbincscu@163.com
\end{document}